\documentclass[oneside, 12pt]{amsart}
\usepackage{amscd, amssymb, amsmath, mathrsfs}
\usepackage[english]{babel}
\usepackage{booktabs}
\usepackage{tikz-cd}
\usepackage[pdftex, colorlinks=true,  citecolor=blue, linkcolor=blue, linktocpage=false]{hyperref}

\newcommand\mylabel[1]{\label{#1}\marginpar{\vspace{-1ex}\medskip\medskip\footnotesize \tt #1}}
\renewcommand\mylabel[1]{\label{#1}}
\newcommand{\mydate}{
\number\day\space
\ifcase\month \or January\or February\or March\or April\or May\or June\or July\or August\or September\or October\or November\or December\fi 
\space\number\year}

\newtheorem{theorem}{Theorem}[section]
\newtheorem*{maintheorem}{Theorem}
\newtheorem{lemma}[theorem]{Lemma}
\newtheorem{proposition}[theorem]{Proposition}
\newtheorem{corollary}[theorem]{Corollary}

\theoremstyle{definition}

\newtheorem*{acknowledgement}{Acknowledgement}

\theoremstyle{remark}

\newcommand{\ZZ}{\mathbb{Z}}

\newcommand{\RR}{\mathbb{R}}

\newcommand{\FF}{\mathbb{F}}

\newcommand{\PP}{\mathbb{P}}
\renewcommand{\AA}{\mathbb{A}}

\newcommand{\ideala}{\mathfrak{a}}

\newcommand{\shN}{\mathscr{N}}
\newcommand{\shL}{\mathscr{L}}

\newcommand{\Frac}{\operatorname{Frac}}

\newcommand{\Kernel}{\operatorname{Ker}}

\newcommand{\lra}{\longrightarrow}

\newcommand{\maxid}{\mathfrak{m}}

\newcommand{\Num}{\operatorname{Num}}

\renewcommand{\O}{\mathscr{O}}

\newcommand{\pdeg}{\operatorname{pdeg}}
\newcommand{\perf}{{\text{\rm perf}}}
\newcommand{\Pic}{\operatorname{Pic}}

\newcommand{\pr}{\operatorname{pr}}
\newcommand{\Proj}{\operatorname{Proj}}

\newcommand{\quadand}{\quad\text{and}\quad}

\newcommand{\ra}{\rightarrow}

\newcommand{\red}{{\operatorname{red}}}
\newcommand{\Reg}{\operatorname{Reg}}

\newcommand{\Sing}{\operatorname{Sing}}

\newcommand{\trdeg}{\operatorname{trdeg}}

\begin{document}

\title[Unirationality]
      {Unirationality and geometric unirationality for hypersurfaces in positive characteristics}

\author[Keiji Oguiso]{Keiji Oguiso}
\address{ Department of Mathematical Sciences, the University of Tokyo, Meguro Komaba 3-8-1, Tokyo, Japan, and National Center for Theoretical Sciences, Mathematics Division, National Taiwan University
Taipei, Taiwan}
\curraddr{}
\email{oguiso@ms.u-tokyo.ac.jp}

\author[Stefan Schr\"oer]{Stefan Schr\"oer}
\address{Mathematisches Institut, Heinrich-Heine-Universit\"at,
40204 D\"usseldorf, Germany}
\curraddr{}
\email{schroeer@math.uni-duesseldorf.de}

\thanks{The first named author is supported by JSPS Grant-in-Aid (S) 15H05738, JSPS Grant-in-Aid (B) 15H03611, and by NCTS Scholar Program.}  
\subjclass[2010]{14M20, 14J70, 14G05, 14J26, 14G17}

\dedicatory{19.\ February 2020}

\begin{abstract}
Building on   work of Segre and  Koll\'ar on   cubic hypersurfaces,
we construct over imperfect fields of characteristic  $p\geq 3$  particular hypersurfaces of degree $p$,
which  show that   geometrically rational schemes that are regular and whose rational points are Zariski dense
are not necessarily unirational. A likewise behaviour holds for 
certain cubic  surfaces in characteristic $p=2$.
\end{abstract}

\maketitle
\tableofcontents

\newcommand{\dashra}{\dashrightarrow}

\section*{Introduction}
\mylabel{Introduction}

Let $F$ be a ground field of arbitrary characteristic $p\geq 0$,
and $X$ be a geometrically integral scheme of dimension $n\geq 0$. One says that $X$ is \emph{rational} or \emph{unirational} if there
is a rational map $\PP^n\dashra X$ that is birational or dominant, respectively.
If this condition holds  after base-change with respect to some finite field extension $F\subset E$, 
one says that $X$ is \emph{geometrically rational} or \emph{geometrically unirational}. 

Let $X\subset \PP^{n+1}$ be an integral cubic hypersurface of dimension $n\geq 2$ that is not a cone.
Generalizing earlier results of Segre \cite{Segre 1943}, Manin \cite{Manin 1986} and Colliot-Th\'el\`ene,
Sansuc and Swinnerton-Dyer \cite{Colliot-Thelene; Sansuc; Swinnerton-Dyer 1987},
Koll\'ar showed over perfect fields $F$ that the following three conditions are equivalent \cite{Kollar 2002}:
\begin{enumerate}
\item The scheme $X$ is unirational.
\item The set of rational points $X(F)$ is non-empty.
\item There is a rational point $a\in X$ whose local ring $\O_{X,a}$ is regular.
\end{enumerate}
For smooth cubic hypersurfaces $X\subset\PP^{n+1}$ , this actually holds 
over arbitrary ground fields $F$.
Furthermore, the   result carries over to imperfect fields of characteristic $p\geq 5$,
and it is asserted that the same holds for the remaining primes under certain technical conditions.

Indeed, Koll\'ar gave the explicit equation $y^3-yz^2+\sum t_ix_i^3=0$
over the function field $F=k(t_1,\ldots,t_n)$ in characteristic three,
which yields a   cubic hypersurface that is regular, geometrically rational and contains exactly three rational points, and
is thus not unirational. He  asks whether a similar equation exists for characteristic two,
and raises for geometrically unirational schemes $X$ the question in what situations the implications
$$
\text{$X$ is unirational}\quad\Longrightarrow\quad 
\text{$X(F)$ is Zariski dense}\quad\Longrightarrow\quad 
\text{$X(F)$ is non-empty} 
$$
might admit reverse implications, say with $X$ smooth and $F$ infinite. 

The goal of this paper is to analyze certain   hypersurfaces $X\subset\PP^{n+1}$ of degree $p$
over imperfect fields $F$
that show that \emph{none of these reverse implications hold, at least with   $X$ regular}.
Generalizing Koll\'ar's equation to arbitrary  $p\geq 3$, we study
$$
y^p-yz^{p-1} +  \sum_{i=1}^nt_ix_i^p=0,
$$
where $x_1,\ldots,x_n,y,z$ are indeterminates and $t_1,\ldots,t_n\in F$ are scalars, with $n\geq 1$.
Here our main result is:

\begin{maintheorem}
{\rm (see Thm.\ \ref{properties hypersurface})}
Suppose the scalars  $t_1,\ldots,t_n\in F$ are algebraically independent over some subfield $k$   of characteristic $p \ge 3$,
and that $F$ is    separable   over the rational function field $k(t_1,\ldots,t_n)$.
Let $F\subset E$ be the extension obtained by adjoining the roots $t_1^{1/p},\ldots,t_n^{1/p}$.
Then our hypersurface $X\subset\PP^{n+1}$    has the following properties:
\begin{enumerate}
\item The scheme $X$ is regular.
\item There is no dominant rational map $\PP^n\dashrightarrow X$   over $F$.
\item The base-change $ X\otimes_FE$ is birational to $\PP^n\otimes_F E$.
\item The set of rational points $X(F)$ is non-empty.
\item If the field $F$ is separably closed, the rational points are  Zariski dense.
\item If $F$ is contained in the field $k((t_1,\ldots,t_n))$,   then   $X(F)$ is finite.
\end{enumerate}
\end{maintheorem}

Properties (i) and (ii) already  hold if the differentials $dt_1,\ldots,dt_n$ in the $F$-vector space of
absolute K\"ahler differentials $\Omega^1_F$
are linearly independent, in other words, if the scalars $t_1,\ldots,t_n\in F$ are \emph{$p$-independent},
a notion going  back to  Teichm\"uller \cite{Teichmueller 1936}. Apparently, this is the correct framework to treat questions
of regularity and unirationality over imperfect fields.

In characteristic $p=2$, we consider the cubic surface $X\subset\PP^3$ defined by the equation
$$
y_1^3+t_1x_1^2y_1 + y_2^3 +  t_2x_2^2y_2 =0
$$
and obtain in Theorem \ref{properties cubic} analogous  results. Here the set of rational points $X(F)$ is always infinite,
because the cubic surface contains a line, but we could not determine whether or not $X(F)$ is  Zariski dense.
As remarked after Proposition \ref{strange properties}, 
this cubic  surface also shows that, for regular cubic hypersurfaces over of characteristic two,  the  implication 
$$
\begin{gathered}
\text{$\exists a\in X(F)$ with   $\O_{X,a}$ regular}\\[-.1cm] 
\text{and   $\pi_a:X\dashrightarrow \PP^n$ separable}
\end{gathered}
\quad\Longrightarrow\quad \text{$X$ is unirational}
$$
formulated in the remark on imperfect ground fields in  \cite{Kollar 2002}, page 468,  does not hold without an  additional assumption. 
The problem seems to be that all tangent plane intersections $C_a=X\cap T_a(X)$, which are non-regular cubic curves,  are 
actually non-integral.
Note that for rational points $a\in X$, the local
ring $\O_{X,a}$ is regular if and only if the scheme $X$ is smooth at the point.

The non-unirationality of our cubic surface depends on the following criterion, which is of independent interest:

\begin{maintheorem}
{\rm (see Thm.\ \ref{frobenius base-change})}
Let $X$ be unirational over some infinite ground field $F$ of characteristic $p>0$. Suppose there  a fibration $f:X\ra\PP^1$ such that the fibers over almost all
rational points $a\in\PP^1$ contain no rational curve. Then the reduced base-change along the relative Frobenius map
$\PP^1\ra\PP^1$ remains unirational.  
\end{maintheorem}

\medskip
The paper is organized as follows:
In Section \ref{Generalities}  we recall basic  facts on $p$-indepen\-dence of scalars $t_1,\ldots,t_n\in F$,
and discuss some implications concerning regularity of schemes and Zariski density of rational points.
In Section \ref{Hypersurfaces} we study hypersurfaces $X\subset\PP^{n+1}$ defined by 
 the equation $y^p-yz^{p-1} +  \sum t_ix_i^p=0$ at odd primes. In Section \ref{Frobenius base-change} we relate unirationality with Frobenius base-change.
This is used in Section \ref{Cubic surface} for the analysis   of  the cubic surface $X\subset\PP^3$ defined by 
the equation  $x_1^3+t_1x_1y_1^2+x_2^3+t_2x_2y_2^2=0$ in characteristic two.

\begin{acknowledgement}
This research was starting during  our stay at the \emph{National Center for Theoretical Sciences}, Taipei.
We would like to thank  Jungkai Chen for the invitation and the NCTS for its hospitality.
The research was also conducted in the framework of the   research training group
\emph{GRK 2240: Algebro-geometric Methods in Algebra, Arithmetic and Topology}, which is funded
by the Deutsche Forschungsgemeinschaft. We also like to thank J\'anos Koll\'ar for helpful comments.
\end{acknowledgement}

\section{Generalities}
\mylabel{Generalities}

Here we recall some  general facts that will be used throughout, concerning K\"ahler differentials, $p$-independence,
regularity, and Zariski density of rational points.
Let $F$ be a    field of characteristic $p>0$, and   
$\Omega^1_{F}=\Omega^1_{F/\ZZ}=\Omega^1_{F/F^p}$
be the $F$-vector space of absolute K\"ahler differentials. The scalars $t\in F$ yield differentials $dt\in\Omega^1_F$,
which form a generating set.
Let us say that a family of scalars $t_i\in F$, $ i\in I$ is \emph{$p$-independent}
if the vectors $dt_i\in \Omega^1_F$ are linearly independent.
We need the following facts:

\begin{proposition}
\mylabel{p-independent}
Consider the following conditions:
\begin{enumerate}
\item The $t_i\in F$   form a separable transcendence basis over a subfield $k$.
\item
The $t_i\in F$ are $p$-independent.
\item  
The $t_i\in F$ are linearly independent over the subfield $F^p$. 
\end{enumerate}
Then  the implications \textrm{(i)$\Rightarrow$(ii)$\Rightarrow$(iii)} hold. Moreover, for each $t\in F$
the condition $dt=0$ is equivalent to $t\in F^p$.
\end{proposition}

\proof
The first implication follows from  \cite{MacLane 1939}, Lemma 3 on page 382.
The second is a consequence of the   characterization of $p$-independence (\cite{Matsumura 1980}, Theorem 86   or \cite{Matsumura 1986}, Theorem 26.5), 
which is frequently taken as a definition: 
The monomials $\prod_{i\in I}t_i^{d_i}\in F$ are linearly independent over the subfield $F^p$, 
where the exponents satisfy $0\leq d_i\leq p-1$
and almost all vanish. In particular, the $t_i\in F$ are linearly independent.

Clearly, each $t\in F^p$ has $dt=0$. Conversely, suppose that $t\in F$ is not a $p$-th power.
The extension $F^p\subset F$ is purely inseparable of height one, so the minimal polynomial of $t$
must be of the form $T^p-\lambda$ for some $\lambda\in F^p$. In turn, the powers $1,t,\ldots,t^{p-1}\in F$
are linearly independent over the subfield $F^p$, and the above characterization shows $dt\neq 0$.
\qed

\medskip
Let us list several elementary but useful permanence properties for $p$-independent scalars:

\begin{proposition}
\mylabel{separable extensions}
Let $F\subset E$ be a separable extension. If  $t_i\in F$, $i\in I$ are $p$-independent,
so are the $t_i\in E$.
\end{proposition}

\proof
According to \cite{Matsumura 1980}, Theorem 88   or \cite{Matsumura 1986}, Theorem 26.6,  the canonical map $\Omega^1_F\otimes_FE\ra\Omega^1_E$ given by
$dt\otimes\lambda\mapsto \lambda dt$  is injective. 
It follows that $F$-linearly independent subsets are mapped to $E$-linearly independent subsets.
\qed

\begin{proposition}
\mylabel{exchange p-independent}
If $t_1,\ldots,t_n\in F$ are $p$-independent, then the same holds for the 
 $t_1,\ldots,t_{n-1},t_n'\in F$ with the new element $t_n'=t_n/t_{n-1}$.
\end{proposition}

\proof
First note that all scalars $t_i$ are non-zero. Set $f=t_{n-1}$ and $g=t_n$. Inside the vector space $\Omega^1_F$, 
the product rule gives
 $g^2d(f/g)=  gdf - f dg$, 
and the assertion follows from the exchange property for linear independent sets.
\qed

\begin{proposition}
\mylabel{remain p-independent}
Suppose that $t_1,\ldots,t_n\in F$ are $p$-independent. Then  the purely inseparable extension
$E=F(t_n^{1/p})$ has  degree $p$,
and the  $t_1,\ldots, t_{n-1}\in E$ remain $p$-independent.
\end{proposition}

\proof
We have $t_n\not\in F^p$, and whence $[E:F]=p$. Clearly, the monomials $t_n^{j/p}$, $0\leq j\leq p-1$
are linearly independent over the subfield  $F$, hence also over $E^p\subset F$, and we infer  that
$\Omega^1_{E/F}$ is one-dimensional, with basis   $dt_n^{1/p}$. 
The field extensions $F^p\subset F\subset E$ gives  an exact sequence
\begin{equation}
\label{upsilon sequence}
0\lra\Upsilon_{E/F/F^p} \lra \Omega^1_{F}\otimes_FE\lra \Omega^1_{E}\lra \Omega^1_{E/F}\lra 0
\end{equation}
Here the term on the left is called the \emph{module of imperfection}, and is defined by the above 
exact sequence; here we follow the notation from \cite{EGA IVa}, Definition 20.6.1.  
Cartier's Equality (\cite{Matsumura 1980}, Theorem 92  or \cite{Matsumura 1986}, Theorem 26.10) 
$$
\dim_E(\Omega^1_{E/F}) = \trdeg_F(E) + \dim_E(\Upsilon_{E/F/F^p})
$$
for the finitely generated field extension $F\subset E$ shows  that our module of imperfection is one-dimensional.
The non-zero vector $dt_n\otimes 1$ clearly belongs to the kernel, whence can be regarded as 
a basis for $\Upsilon_{E/F/F^p}$. 
It follows that the remaining vectors $dt_1,\ldots,dt_{n-1}$ remain linearly independent in $\Omega^1_E$.
\qed

\medskip
Now let $x_0,\ldots,x_n$ be indeterminates for some $n\geq 0$, and regard $\PP^n$ as the homogeneous spectrum
of the polynomial ring $F[x_0,\ldots,x_n]$. Given a sequence of scalars $t_0,\ldots,t_n\in F$, not all of which vanish,
we consider the Fermat hypersurface $D\subset \PP^n$  
defined by the equation $t_0x_0^p+\ldots+t_nx_n^p=0$.
Note that $D$ is irreducible but  geometrically non-reduced, and becomes a $p$-fold hyperplane
after base-changing to the perfect closure.  

\begin{proposition}
\mylabel{fermat hypersurface}
Suppose $t_0=1$. Then the  scheme $D$ is regular if and only if the $t_1,\ldots,t_n\in F$ are $p$-independent.
\end{proposition}

\proof
The extension $F'=F^p(t_1,\ldots,t_n)$ defines  an intermediate field   $F^p\subset F'\subset F$.  
The \emph{$p$-degree} 
$
d=\pdeg(F'/F^p)
$
is defined as the vector space dimension of $\Omega^1_{F'/F^p}$,
and is also characterized by the degree formula $[F':F^p]=p^d$.
We have $d\leq n$, because  the differentials $dt_1,\ldots,dt_n\in\Omega^1_{F'/F^p}$ form a generating set.
According to \cite{Schroeer 2010}, Theorem 3.3 the scheme $D$ is regular if and only if $d=n$.
Hence we have to show the equality
\begin{equation}
\label{vector space dimensions}
\dim_{F'}(\Omega^1_{F'/F^p})=\dim_F( Fdt_1+\ldots +Fdt_n)
\end{equation}
of vector space dimensions. Taking $p$-th roots, we see that the left hand side equals the dimension
of $ \Omega^1_{E/F} $. Here    $F\subset E$ denotes the extension generated by $t_1^{1/p},\ldots, t_n^{1/p}$, to avoid confusion with $F'$.
Using induction on $n\geq 0$ with Proposition \ref{remain p-independent}, 
one sees that the right hand side $r=\dim_F( Fdt_1+\ldots +Fdt_n)$
obeys the formula $[E:F]=p^r$, hence also 
coincides with the dimension of $ \Omega^1_{E/F}$. This gives the desired equality \eqref{vector space dimensions}.
\qed

\medskip
Now suppose that $X$ is an $F$-scheme of finite type.
One says that $X$ is \emph{geometrically reduced} if for some  algebraically  closed field extension $E$,
the base-change $X'=X\otimes_FE$ is reduced.

\begin{lemma}
\mylabel{zariski dense}
If the scheme $X$ is geometrically reduced and the field 
$F$ is separably closed,  then the set of rational points $X(F)$ is Zariski dense.
\end{lemma}

\proof
We have to verify that each non-empty open set contains a rational point,
so it suffices to check that $X(F)$ is non-empty, and we may assume that $X$ is affine.
By Bertini's Theorem (\cite{Jouanolou 1983}, Theorem 6.3) there is a hyperplane $H\subset X$
that remains geometrically reduced.
By induction on the dimension, this reduces us to the case $\dim(X)=0$. Hence our scheme is the spectrum 
of a product $E_1\times\ldots\times E_r$ of $r\geq 1$ separable field extensions.
Since $F$ is separably closed, we must have $E_i=F$.

The following more direct argument was suggested to us by J\'anos Koll\'ar: 
According to \cite{MacLane 1939}, Theorem 15 the function field of $X$ has a a separating transcendence basis over $F$.
In turn, we may assume that $X$ is \'etale over $\AA^n$. For each rational point $a\in\AA^n$ lying in the
image of $X$, the preimage is the spectrum of a product $E_1\times\ldots\times E_r$ as above.
\qed

\medskip 
Suppose now that $X$ is  equidimensional of dimension $n\geq 0$. Then 
the \emph{locus of non-smoothness} $\Sing(X/F)$ is the set of points $a\in X$ where $\Omega^1_{X/F}\otimes \kappa(a)$ has
vector space dimension $d>n$. It has a natural scheme structure, defined via Fitting ideals for
the coherent sheaf $\Omega^1_{X/F}$, compare the discussion in \cite{Fanelli; Schroeer 2020}, Section 2.
Depending on the context, we also call $\Sing(X/F)$ the \emph{scheme of non-smoothness}.

\begin{lemma}
\mylabel{regular}
Suppose that $\Sing(X/F)$ and some effective Cartier divisor $D\subset X$
have the same support, and that $X$ contains no embedded components. Then $X$ is geometrically reduced but geometrically non-normal.
Furthermore, the scheme $X$ is regular provided that $D$ is regular. 
\end{lemma}

\proof
The open set $X\smallsetminus D$ is smooth. The base-change $X'=X\otimes_FE$ to the perfect closure $E=F^\perf$
also contains no embedded component, and is generically smooth. In turn, the structure sheaf $\O_{X'}$ has
no  non-zero nilpotent elements, so $X$ is geometrically reduced.
Let $\zeta $ be some generic point in $D'=D\otimes_FE$. Then the local ring $\O_{X',\zeta}$ is one-dimensional and not regular.
By Serre's Characterization, the scheme $X$ is not geometrically normal.

Suppose now that the scheme $D$ is regular.
Fix a point $a\in D$, and let $f\in\O_{X,a}$ be an element  defining the Cartier divisor
in some neighborhood. This element is regular and contained in the maximal ideal.
Since the local ring $\O_{D,a}=\O_{X,a}/(f)$ is regular, the same must hold for $\O_{X,a}$.
\qed

\section{Hypersurfaces of \texorpdfstring{$p$}{p}-degree}
\mylabel{Hypersurfaces}

Let $F$ be a  ground field of characteristic $p\geq 3$. Fix some integer $n\geq 1$ and   scalars $t_1,\ldots,t_n\in F$,
only subject to the condition $t_1\neq 0$. 
Regard $\PP^{n+1}$ as the homogeneous spectrum of the polynomial ring
$F[x_1,\ldots,x_n,y,z]$.
We now consider the hypersurface $X\subset\PP^{n+1}$ of dimension $\dim(X)=n$ and degree $\deg(X)=p$ defined by the equation
\begin{equation}
\label{kollar equation}
y^p-yz^{p-1} + \sum_{i=1}^n t_ix_i^p=0.
\end{equation}
For function fields $F=k(t_1,\ldots,t_n)$ in characteristic three, 
this is the cubic hypersurface studied by Koll\'ar in \cite{Kollar 2002}, Section 4.
Here we work over arbitrary characteristics $p\geq 3$ and more general ground fields $F$.

\begin{proposition}
\mylabel{geometrically integral}
The scheme $X$ is geometrically integral.
\end{proposition}

\proof
Replacing $F$ by some algebraic closure, we have to show that the left-hand side of \eqref{kollar equation} is an 
irreducible polynomial. Set $x=\sum t_i^{1/p}x_i$ and $v=x+y$. 
Now our task is to verify that $P(v)=v^p-yz^{p-1}$ is irreducible as polynomial over $R=k[y,z]$.
This follows immediately with the Eisenstein Criterion with the prime element $y\in R$.
\qed

\begin{proposition}
\mylabel{birational}
If $t_1,\ldots,t_n\in F^p$, then the scheme $X$ is birational to $\PP^n$.
\end{proposition}

\proof
As in the previous proof, we may assume that our hypersurface $X\subset\PP^{n+1}$ is given by
the equation $y^p-yz^{p-1} + x_1^p=0$. This does not involve the variables $x_2,\ldots,x_n$,
hence $X$ is a cone with respect to the $(n-2)$-dimensional linear subspace $V\subset\PP^{n+1}$ given by $x_1=y=z=0$ as apex,
over the plane curve $C\subset\PP^2$ defined by the equation $x^p-yz^{p-1}=0$, where we have made the substitution $x=y+x_1$.

Geometrically, this means that $X$ is birational to $C\times\PP^{n-1}$, and it remains to check that the integral curve $C$ is rational.
On the affine chart given by $z\neq 0$, the coordinate ring for the curve becomes
the polynomial ring $F[x/z]$, hence  $C$ must be rational.
\qed

\begin{proposition}
\mylabel{regular hypersurface}
The scheme of non-smoothness $\Sing(X/F)\subset X$ and the effective Cartier divisor $D\subset X$
defined by the equation $z=0$ have the same support. 
Moreover,   $X$ is regular  provided that  $t_1,\ldots,t_n\in F$ are   $p$-independent.
\end{proposition}

\proof
For our hypersurface $X\subset\PP^{n+1}$, the scheme of non-smoothness $\Sing(X/F)$ is defined by
the additional equations coming from the partial derivatives of \eqref{kollar equation}.
These partial derivatives are $z^{p-1}$ and $-yz^{p-2}$. It follows that $D$ and $\Sing(X/F)$
have the same support.

Now suppose that $t_1,\ldots,t_n\in F$ are $p$-independent. We may regard $D$ as the divisor
in $\PP^n$ defined by the Fermat equation $y^p +t_1x^p+\ldots+t_nx_n^p$.
According to Proposition \ref{fermat hypersurface}, the hypersurface $D$ is regular.
By Lemma \ref{regular}, the scheme $X$ is regular as well.
\qed

\medskip
In order to apply induction, we will relate our hypersurface in dimension $n$ with one in dimension $n-1$.
This is based on the following observation:

\begin{lemma}
\mylabel{projective equivalence}
Suppose  $n\geq 2$, that $t_{n-1}\neq 0$ and that $t_n/t_{n-1}\in F^p$. Then the hypersurface $X\subset\PP^{n+1}$
is projectively equivalent to the hypersurface $X'\subset\PP^{n+1}$ defined by another equation of the form \eqref{kollar equation},
with  coefficients $t_i'=t_i$ for $i\leq n-1$ and $t'_n=0$.
\end{lemma}

\proof
Let $\lambda\in F$ be the scalar with $\lambda^p=t_n/t_{n-1}$, rewrite the   equation \eqref{kollar equation} as
$$
y^p-yz^{p-1} + t_1x_1^p+\ldots+t_{n-2}x_{n-2}^p + t_{n-1}(x_{n-1}  + \lambda  x_n )^p=0,
$$
and use the coordinate change $x_{n-1}'=x_{n-1}  + \lambda  x_n $.
\qed

\begin{proposition}
\mylabel{not unirational}
If $t_1,\ldots,t_n\in F$ are $p$-independent, then the scheme $X$ is not unirational.
\end{proposition}

\proof
We proceed by induction on $n=\dim(X)$.
Suppose first that $n=1$. Seeking a contradiction, we assume that there is a rational dominant map $\PP^1\dashrightarrow X$.
By L\"uroth's Theorem, $X$ is  birational to $\PP^1$.  According to Proposition \ref{regular hypersurface}, the curve $X$ is regular,
so we actually have an isomorphism $X\simeq\PP^1$, in particular   $X$ is smooth. On the other hand, the scheme
of non-smoothness $\Sing(X/F)$ is non-empty, contradiction.

Suppose now that $n\geq 2$, and that the assertion is true for $n-1$.
Seeking a contradiction, we assume that there is a rational dominant map $\PP^n\dashrightarrow X$.
Let us write $X=X_F(t_1,\ldots,t_n)$ to indicate the dependence of our hypersurface $X\subset\PP^n$
on the ground field $F$ and the scalars $t_1,\ldots,t_n\in F$.
Consider its base-change $\PP^n_E\dashrightarrow X_E(t_1,\ldots,t_n)$ for the field extension $E=F(t_n^{1/p})$.
According to  Lemma \ref{projective equivalence} there is linear isomorphism $X_E(t_1,\ldots,t_n)\ra X_E(t_1,\ldots,t_{n-1},0)$.
The latter becomes a cone in $\PP^{n+1}_E$, because its  equation no longer involves the indeterminate $x_n$,  whence there is a 
dominant rational map 
$$
X\otimes_FE=X_E(t_1,\ldots,t_{n-1},0)\dashrightarrow X_E(t_1,\ldots,t_{n-1})=X'.
$$
Composing these maps we get  a dominant rational map $\PP^n_E\dashrightarrow X'$.
According to  \cite{Kollar 2002}, Lemma 2.3 the hypersurface  $Y$ is unirational.
On the other hand, the scalars $t_1,\ldots,t_{n-1}\in E$ are $p$-independent according to Proposition \ref{remain p-independent}.
By induction hypothesis, the hypersurface $X'\subset\PP^n_E$ is not unirational, contradiction.
\qed

\medskip
The hypersurface $X\subset \PP^{n+1}$ contains the  obvious rational points  
\begin{equation}
\label{obvious points}
(0:\ldots:0:\lambda:1), \quad \lambda\in\FF_p.
\end{equation}  
Under suitable assumptions on the ground field $F$, there are no further rational points:

\begin{proposition}
\mylabel{rational points}
Suppose that $F$ is contained in the field $k((t_1,\ldots,t_n))$  of formal Laurent series  
with respect to indeterminates $t_1,\ldots,t_n$ and some subfield $k$. Then  $X(F)$ consists of the  $p$ rational  points listed in \eqref{obvious points}.
\end{proposition}

\proof
This is essentially Koll\'ar's argument from \cite{Kollar 2002}, Section 4, which we repeat for the convenience of the reader.
It suffices to treat the case that $F$ equals the field of formal Laurent series over an infinite field $k$.
This means $F=\Frac(R)$ for the ring   $R=k[[t_1,\ldots,t_n]]$.
Let  $a\in X(F)$ be a rational point, and write it  as  $a=(h_1:\ldots:h_n:f:g)$ with some
relatively prime power series $h_i,f,g\in R$. This is indeed possible  because the ring $R$ is factorial by  \cite{Matsumura 1986}, Theorem 20.8.
Our task is to show that the $h_i$ vanish. Seeking a contradiction, we assume that this is not the case. 
Given some exponents $u_i\geq 1$, we obtain 
a homomorphism $\varphi: R\ra k[[t]]$ defined by   $t_i\mapsto t^{u_i}$, inducing 
an equation $f^p-fg^{p-1}+th^p=0$, now with $f,g,h\in k[[t]]$.
According to \cite{AC 7}, \S3, No.\ 7, Lemma 2 we may choose the exponents so that $h\neq 0$. Then also $f\neq 0$.

Dividing by some common factor, we may assume that $\gcd(f,g,h)=1$. Each irreducible factor $d$ of $\gcd(f,g)$
has the property $d^p| th^p$. Since $t$ is a prime element, we must have  $d|h$, contradiction. Thus $\gcd(f,g)=1$.
Rewrite our equation as
$
th^p=\prod_{j=0}^{p-1}(f-jg) 
$.
The factors $P_j=f-jg$ on the right are pairwise coprime, because this holds for $f,g$.
Hence we can write $f-jg=Q_j^p$ for all $j$  with   one exception $i$, which has $f-ig = tQ^p_i$.
Then 
$$
tQ_i^p + (\sum_{j\neq i}Q_j)^p  = \sum_{j=0}^{p-1}(f-jg) = pf - p\frac{p-1}{2} g = 0.
$$
We conclude that in the prime factorization of $tQ_i^p$, all exponents are divisible by $p$.
This contradicts the fact that  $t$ is a prime element in the ring $k[[t]]$.
\qed

\medskip
We now summarize our results in the following form:
 
\begin{theorem}
\mylabel{properties hypersurface}
Suppose the  scalars $t_1,\ldots,t_n\in F$ are algebraically independent over some subfield $k$  of characteristic $p \ge 3$,
and that $F$ is    separable   over the rational function field $k(t_1,\ldots,t_n)$.
Let $F\subset E$ be the extension obtained by adjoining the roots   $t_1^{1/p},\ldots,t_n^{1/p}$. 
Then the hypersurface $X\subset\PP^{n+1}$  that is defined by the equation $y^p-yz^{p-1} + \sum_{i=1}^n t_ix_i^p=0$  has the following properties:
\begin{enumerate}
\item The scheme $X$ is regular.
\item There is no dominant rational map $\PP^n\dashrightarrow X$  over $F$. 
\item The base-change $ X\otimes_FE$ is birational to $\PP^n\otimes_F E$.
\item The set of rational points $X(F)$ is non-empty.
\item If the field $F$ is separably closed, the rational points are  Zariski dense.
\item If $F$ is contained in the field $k((t_1,\ldots,t_n))$,   then   $X(F)$ is finite.
\end{enumerate}
\end{theorem}

\proof
According to Proposition \ref{p-independent}, the scalars $t_1,\ldots,t_n\in F$ are $p$-independent,
so the scheme $X$ must by regular by Proposition \ref{regular hypersurface}.
Furthermore, it is not unirational according to Proposition \ref{not unirational}. 
The base-change $X\otimes_FE$ becomes rational, in light of Proposition \ref{birational}.
If $F$ is separably closed, the rational points  must be dense by Lemma \ref{zariski dense}.
If  $F$ is contained in the field of formal Laurent series, we saw in Proposition \ref{rational points} that there
are only $p$ rational points.
\qed

\medskip
With the setting of  the above theorem, our regular scheme  $X$ is  geometrically unirational
but not unirational. 
Furthermore,  no separable extension achieves unirationality.
Moreover, we see that the reverse implications
$$
\text{$X$ is unirational}\quad\Longrightarrow\quad 
\text{$X(F)$ is Zariski dense}\quad\Longrightarrow\quad 
\text{$X(F)$ is non-empty} 
$$
for geometrically unirational schemes $X$ over infinite fields $F$
discussed by Koll\'ar in \cite{Kollar 2002}, Question 1.3 do not hold for regular schemes.

\section{Unirationality and Frobenius base-change}
\mylabel{Frobenius base-change}

Let $F$ be an infinite ground field of characteristic $p>0$.
Suppose $X$ is an integral proper scheme endowed with a surjective morphism $f:X\ra \PP^1$.
Write the projective line as the homogeneous spectrum of $F[T_0,T_1]$, and regard the indeterminates
$T_i$ as global sections of the ample sheaf $\O_{\PP^1}(1)$.
Fix an integer $\nu\geq 1$. The resulting  global sections $T_i^{p^\nu}$ of $\O_{\PP^1}(p^\nu)$ define
a purely inseparable morphism $h:\PP^1\ra\PP^1$ of degree $p^\nu$. This map can also be described
by the inclusion of coordinate rings   $F[s^{p^{\nu}}]\subset F[s]$, where we set $s=T_1/T_0$.  
This reveals that $h:\PP^1\ra\PP^1$ coincides with the iterated relative Frobenius map for the projective line.
Let us write $X'=(X\times_{\PP^1}\PP^1)_\red$ for the ensuing base-change, endowed with the reduced scheme structure.

\begin{theorem}
\mylabel{frobenius base-change}
Suppose the scheme $X$ is unirational, and that for almost every rational point $a\in\PP^1$, the fiber $f^{-1}(a)$ contains no rational curve. 
Then the reduced base-change $X'=(X\times_{\PP^1}\PP^1)_\red$ is unirational as well.
\end{theorem}

\proof
Set $n=\dim(X)$, and choose a dominant rational map $\PP^1\times\PP^{n-1}\dashrightarrow X$.
By the Valuative Criterion for properness, the domain of definition contains   $\PP^1_U=\PP^1\times U$
for some open dense set $U\subset\PP^{n-1}$, so we have a dominant morphism $g:\PP^1_U\ra X$.

We now write $B=\PP^1$ for the base of the given surjection $f:X\ra\PP^1=B$. Let $b_1,\ldots,b_r\in B$
be the finitely many rational points whose fibers contain rational curves.
The preimages of $f^{-1}(b_i)$ on $\PP^1_U$ are closed sets not containing the generic point.
Since the projection $\PP^1_U\ra U$ is proper, we may shrink $U$ and suppose that the image of $g:\PP^1_U\ra X$ is
disjoint from the fibers $f^{-1}(b_i)$. This means that for every  rational point  $u\in U$, the image
$g(\PP^1_u)\subset X$ is not contained in any of the  fibers of $f:X\ra B$, and thus dominates $B$.
It follows that for the generic point $\eta\in U$, the induced projection $\PP^1_E=\PP^1_\eta\ra B=\PP^1$ is surjective,
where $E=\kappa(\eta)$ denotes the function field of the open set $U\subset\PP^{n-1}$.

Consider the composite morphism $\PP^1_U\ra B$ and the ensuing base-change $(\PP^1_U)\times_B B$
with respect to the purely inseparable morphism $h:B=\PP^1\ra\PP^1=B$ of degree $\deg(h)=p^\nu$.
 It comes with a projection $\pr:(\PP^1_U)\times_BB\ra U$ and a dominant morphism $(\PP^1_U)\times_BB\ra X\times_BB$.
To check that $X'$ is unirational, it thus suffices to verify that the reduction of the generic fiber $\pr^{-1}(\eta)$ 
is a rational curve over the function field $E=\kappa(\eta)$ of the open set $U\subset\PP^{n-1}$.

This is a consequence of the following property  of the iterated relative Frobenius map $h:\PP^1\ra\PP^1$: We claim that
for each field extension $F\subset E$ and each surjective $F$-morphism $\varphi:\PP^1_E\ra\PP^1_F$, there is a commutative diagram
$$
\begin{CD}
\PP^1_E			@<h_E<<	\PP^1_E\\
@V\varphi VV			@VV\psi V\\
\PP^1_F			@<<h<	\PP^1_F
\end{CD}
$$
for some $\psi$. Indeed, the morphism $\varphi$ is defined via some invertible sheaf $\shL=\O_{\PP^1_E}(n)$
and two global sections without common zeros, which can be viewed
as homogeneous polynomials $Q_0,Q_1\in E[T_0,T_1]$ of degree $n$ that are relatively prime. 
Set $q=p^\nu$.
Then the morphism $\psi$   defined by the polynomials $Q_0^q,Q_1^q\in E[T_0^q,T_1^q]$
makes the diagram commutative.

The above diagram yields a surjection $\PP^1_E \ra \PP^1_E\times_{\PP^1_F}\PP^1_F$. This is an $E$-morphism,
because the iterated relative Frobenius map $h_E$ is an $E$-morphism. L\"uroth's Theorem (see for example \cite{Kollar; Smith; Corti 2004}, Proposition 1.10)
ensures that the reduction of the fiber product is a rational curve over $E$.
\qed

\medskip
The following consequence will later play an important role:

\begin{corollary}
\mylabel{criterion not unirational}
Suppose that for almost every rational point $a\in\PP^1$, the fiber $f^{-1}(a)$ contains no rational curves,
and that $X'$ is birational to $Z\times\PP^{n-1}$, where $Z$ is not a rational curve. Then $X$ is not unirational.
\end{corollary}

\proof
Seeking a contradiction, we assume that $X$ is unirational. By the theorem, $X'$ is unirational and hence $Z$
are rational, contradiction.
\qed

\section{A cubic  surface  in characteristic two}
\mylabel{Cubic  surface}

Let $F$ be a ground field of characteristic $p=2$. 
Regard $\PP^3$   as the homogeneous spectrum of the polynomial ring $F[x_1,x_2,y_1,y_2]$,
and let $t_1,t_2\in F$ be scalars, subject only to the condition $t_1\neq 0$ and $t_2 \not= 0$.
The goal of this section is to study the cubic surface $X\subset\PP^3$ defined by the equation
\begin{equation}
\label{cubic equation}
y_1^3+t_1x_1^2y_1 + y_2^3 +  t_2x_2^2y_2 =0.
\end{equation}
The defining polynomial is irreducible, which can be seen by setting $x_2=0$ and observing that
$y_1(y_1^2+t_1x_1^2)$ is not a cube in $F[x_1,y_1]$. Thus our del Pezzo surface $X$ is a geometrically integral.

In what follows we shall see that if the scalars are $p$-independent, the scheme $X$ is regular and geometrically rational, 
yet not unirational. Moreover, it has Picard number $\rho=2$ and contains exactly one $(-1)$-curve $L$.
We do not now whether or not $X(F)$ is Zariski dense.
 
\begin{proposition}
\mylabel{locus non-smoothness}
The scheme of non-smoothness $D=\Sing(X/F)$  is  an irreducible curve defined inside $\PP^3$   by the two equations
$y_1^2+t_1x_1^2=0$ and $y_2^2+  t_2x_2^2=0$. Moreover, the inclusion $D\subset X$ is Cartier.
\end{proposition}
 
\proof
The partial derivatives of the defining polynomial  $P=y_1^3+t_1x_1^2y_1 + y_2^3 +  t_2x_2^2y_2$
with respect to $y_i$ are $P_i = y_i^2+t_ix_i^2$, whereas $\partial P/\partial x_i=0$.
Moreover, the jacobian ideal $\ideala=(P,P_1,P_2)$ is already generated by the two partial derivatives,
which yields the assertion on the embedding $D\subset \PP^3$. 
If $t_i\in F$ are squares, a change of coordinate reveals that $D$ is the intersection 
of two double planes, which shows that $D$ is an irreducible curve.

From \eqref{cubic equation}, one sees that on the open set given by $y_2\neq 0$,
the inclusion $D\subset X$ is already defined by the single equation $y_1^2+t_1x_1^2=0$.
An analogous statement holds on the open set given by $y_1\neq 0$.
It follows that $D\subset X$ is Cartier outside the closed set $L\subset X$ defined by $y_1=0$ and $y_2=0$.
From the equations for $D\subset\PP^3$ one sees it is disjoint from $L$, hence $D\subset X$ must be Cartier.
\qed
 
\medskip
The line $L\subset\PP^3$ given by the equations $y_1=0$ and $y_2=0$  is contained in $X$  and   lies in the smooth locus.
 The adjunction formula
for the inclusions $X\subset\PP^3$ and $L\subset X$ gives $\omega_X=\O_X(-1)$ and 
$-2=(L+K_X)\cdot L=L^2-1$. Hence:

\begin{proposition}
\mylabel{selfintersection}
The selfintersection number of the line $L$ on the cubic surface $X$ is given by $L^2=-1$.
In other words, $L\subset X$ is a $(-1)$-curve.
\end{proposition}

Now consider the plane  $H_1\subset\PP^3$ given by the equation $y_1=0$.
Then the plane section $H_1\cap X$ is defined by $y_1=0$ and $y_2(y_2^2+t_2x_2^2)=0$, thus decomposes as $L+C_1$,
where $C_1$ is the irreducible conic defined by $y_1=0$ and $y_2^2+t_2x_2^2=0$. Likewise, 
the plane $H_2\subset \PP^3$ defined by $y_2=0$ has  $H_2\cap X = L+C_2$, where
the irreducible conic $C_2$ is defined by $y_2=0$ and   $y_1^2+t_1x_1^2=0$.   

The equations reveal that $C_1\cap C_2=\varnothing$.
Moreover, the curves  $C_i\subset X$ are Cartier, because the intersections $C_i\cap L$ lies in the smooth locus. 
Since $H_1,H_2\subset \PP^2$ are linearly equivalent, the same holds for $C_1,C_2\subset X$. In turn, the invertible sheaf 
$\shL=\O_X(C_1)$ is globally generated, and the two-dimensional linear  system inside $H^0(X,\shL)$ 
generated by global sections defining  $C_i\subset X$ yield a morphism $f:X\ra\PP^1$
with $\shL=f^*\O_{\PP^1}(1)$.

Now it is convenient to  use the term \emph{double line} for a curve isomorphic to 
the first infinitesimal neighborhood of a line $\PP^1$ in $\PP^2$.   Note that the \emph{twisted forms of the double line} are precisely the
conics that are geometrically non-reduced.

\begin{proposition}
\mylabel{fibration}
The morphism $f:X\ra\PP^1$ extends the rational map $X\dashra \PP^1$ given by $(x_1:y_1:x_2:y_2)\mapsto (y_1:y_2)$.
All fibers  are twisted forms of the  double line.
The induced finite morphisms
$$
f:L\lra\PP^1\quadand f:D=\Sing(X/F)\lra\PP^1
$$
are purely inseparable of degree two and four, respectively.
\end{proposition}

\proof
Let $s_1,s_2$  be     sections of $\shL$ defining  $C_1,C_2\subset X$, and $E\subset H^0(X,\shL)$
the resulting linear system.
By construction, we have $\shL = \O_X(1)\otimes \O_X(-L)$. Under the canonical inclusion $\shL\subset\O_X(1)$
and up to scalars,
the sections $s_i$ become the restrictions of $y_i\in H^0(\PP^3,\O_{\PP^3}(1))$, and $L$ is the fixed part
of the $y_1,y_2$.
The rational map $\varphi:\PP^3\dashra\PP^1$ given by $(x_1:y_1:x_2:y_2)\mapsto (y_1:y_2)$
has the open set $U=\PP^3\smallsetminus L$ as domain of definition,  and it also can be described by the two-dimensional linear
system generated by  $y_1,y_2\in H^0(\PP^3,\O_{\PP^3}(1))$. Thus the map $\varphi|X$
coincides with the morphism $f:X\ra \PP^1$ on the open set $X\cap U$.

Now let $a\in\PP^1$ be a   point. To check that the fiber is a twisted form of the double line,
it suffices to treat the case that $a=(\lambda_1:\lambda_2)$ is a rational point.
Then the fiber $Z=f^{-1}(a)$ is the zero-scheme for $\lambda_1s_1+\lambda_2 s_2$, and is contained in the
zero-scheme $Z'\subset X$ for $\lambda_1y_1+\lambda_2y_2$, which is a plane section.
In turn, $Z'=Z\cup L$ is a reducible cubic curve, thus decomposes into the union of a conic $Z$ and a line $L$.
This shows that the  fiber $Z=f^{-1}(a)$ is isomorphic to a conic.
To proceed,  it suffices by symmetry to treat the case that $\lambda_2=1$, and we write $\lambda=\lambda_1$. 
Then $f^{-1}(a)\subset X $ is defined inside $\PP^3$ by the 
homogeneous   equations
\begin{equation}
\label{fiber computation}
\lambda y_1+y_2=0\quadand (1+\lambda^3)y_1^2+t_1x_1^2 +\lambda t_2x_2^2=0,
\end{equation}
which indeed is a twisted form of the double line.
Taking intersections with $L$ and $D=\Sing(X/F)$, one sees that the induced projections are purely inseparable
of degree $d=2$ and $d=4$, respectively.
\qed

\medskip
We now come to the main result on our cubic surface:

\begin{theorem}
\mylabel{properties cubic}
Suppose the scalars $t_1,t_2\in F$ are $p$-independent. Let $F\subset E$ be the field extension obtained by
adjoining the root  $\sqrt{t_1}$. Then the cubic surface $X\subset\PP^3$
defined by the equation $y_1^3+t_1x_1^2y_1 + y_2^3 +  t_2x_2^2y_2 =0$  has the following properties:
\begin{enumerate}
\item The scheme $X$ is regular.
\item There is no   dominant rational map $\PP^2\dashra X$  over $F$. 
\item The base-change $X\otimes_FE$ is birational to $\PP^2\otimes_FE$.
\item The set of rational points $X(F)$ is infinite.
\item If $F$ is separably closed, the rational points are Zariski dense.
\end{enumerate}
\end{theorem}

\proof
The assertion (iv) is a consequence of Proposition \ref{zariski dense}, and (iv) follows from the existence
of the line $L\subset X$.
Over the field extension $E$, we set $x_1'=y_1+\sqrt{t_1}x_1$.
In the new indeterminates $x_1',y_1,x_2,y_2$ our cubic surface is given by the equation
$y_1x_1^2 +  y_2^3 +  t_2x_2^2y_2 =0$. Localizing with respect to $x_1$ we see that $y_1$
can be expressed by the other three indeterminates. This ensures that the base-change $X\otimes_FE$
is a rational surface, hence (iii).

We next verify that the scheme $X$ is regular. Recall that the scheme of non-smoothness $D=\Sing(X/F)$ was described
in Proposition \ref{locus non-smoothness}. Consider first the non-rational closed point
$a=(1:0:\sqrt{t_1}:0)\in   D$.
On the open set given by $x_1\neq 0$, the cubic surface is defined by the inhomogeneous equation
$$
\frac{y_1}{x_1} \left( \left(\frac{y_1}{x_1}\right)^2 + t_1\right) + 
\left(\frac{y_2}{x_1}\right)^3 + t_2\left(\frac{x_2}{x_1}\right)^2  \frac{y_2}{x_1} =0,
$$
and the polynomial on the left lies in the maximal ideal of $\maxid_R$ of the local ring $R=\O_{\AA^3,a}$, 
but not in $\maxid_R^2$. In turn, $\O_{X,a}$ is regular. By symmetry, the same
holds at the closed point $b=(0:1:0:\sqrt{t_2})$.
According to  Lemma  \ref{regular}, it suffices to verify that the scheme $D\smallsetminus\{a,b\}$
is regular.  This lies in the open set given by $y_1,y_2\neq 0$,
hence equals the spectrum of the ring
$$
F[u,v,w^{\pm 1}]/(1+t_1u_1^2,1+t_2u_2^2)
$$
where we set $u_1=x_1/y_1$ and $u_2=x_2/y_2$ and $w=y_1/y_2$. Clearly, this ring is isomorphic
to the ring of Laurent polynomials in $w$ over the tensor product $A=F(\sqrt{t_1})\otimes_FF(\sqrt{t_2})$.
The latter is a field, because $t_1,t_2\in F$ are $p$-independent, hence $D\smallsetminus\{a,b\}$ is indeed regular.
This establishes (i).

It remains to verify   (ii), which is the most interesting part. 
For this we apply Corollary \ref{criterion not unirational} to our  fibration $f:X\ra\PP^1$.
Let us examine the fiber $f^{-1}(a)$ over the rational points $a=(\lambda:1)$ with $\lambda^3\neq 1$,
which means $a\not\in\PP^1(\FF_4)$.
According to \eqref{fiber computation} this is a conic $C\subset\PP^2_F$ given by the equation
\begin{equation}
\label{fiber conic}
(1+\lambda^3)u_0^2 + t_1u_1^2 + \lambda t_2u_2^2 =0
\end{equation}
in some indeterminates $u_0,u_1,u_2$. Base-changing to the field extension $F'=F(\sqrt{\lambda})$,
and making a linear change of variables, the equation can be rewritten as 
\begin{equation}
\label{modified conic}
v_0^2 + t_1v_1^2 + t_2v_2^2=0.
\end{equation}
The short exact sequence \eqref{upsilon sequence} and Cartier's Equality (\cite{Matsumura 1980}, Theorem 92   or \cite{Matsumura 1986}, Theorem 26.10)
reveal that  the kernel for  $\Omega^1_{F}\otimes F'\ra\Omega^1_{F'}$ is at most one-dimensional.
So without restriction, we may assume that $dt_1\in\Omega^1_{F'}$ remains non-zero.
According to \cite{Schroeer 2010}, Theorem 3.3 the conic $C\otimes_FF'$ is reduced, hence the same
holds for $C$. Since the latter is geometrically non-reduced, it is not rational. 
Summing up, for almost all rational points $a\in \PP^1$, the fiber $f^{-1}(a)$ is not rational.

We proceed with a similar computation for the  generic fiber  of $f:X\ra\PP^1$ and its Frobenius base-change.
Regard $\PP^1$ as the homogeneous spectrum
of  $F[y_1,y_2]$,  and now write $\lambda=y_2/y_1$ for the transcendental generator of the function field.
Then the generic fiber for $f:X\ra \PP^1$ is the conic 
given by \eqref{fiber conic} over $F(\lambda)$, and the generic fiber of the Frobenius base-change is given by
the same equation over $F(\sqrt{\lambda})$.
This is already defined over the subfield $F$, and we conclude that the Frobenius base-change
$X\times_{\PP^1}\PP^1$ is birational to $C\times\PP^1$,
where $C\subset\PP^2_F$ is the conic defined by the above equation.
According to Proposition \ref{fermat hypersurface}, the curve $C$ is regular.
Being geometrically non-reduced, it is not rational. Thus Corollary \ref{criterion not unirational} applies, and
we conclude that $X$ is not unirational.
\qed

\medskip
Each rational point $a\in X\subset\PP^3$ comes from a linear surjection $\varphi:F^4\ra F$. 
Choosing  $\Kernel(\varphi)\simeq F^3$ we obtain  a rational map $\pi_a:X\dashrightarrow \PP^2$.
If moreover $\O_{X,a}$ is regular,  
the  intersection $C_a=X\cap T_a(X)$ is a singular cubic curve in the tangent plane $T_a(X)\subset\PP^3$.
Note that these $C_a\subset T_a(X)$ are crucial in the work of Segre \cite{Segre 1943}, Manin \cite{Manin 1986}
and Koll\'ar \cite{Kollar 2002}.
 
\begin{proposition}
\mylabel{strange properties}
The rational map $\pi_a:X\dashra\PP^2$ is purely inseparable if and only if $a\in L$. Moreover, the intersection $C_a$
is not integral for every $a\in\Reg(X)$.
\end{proposition}
 
\proof
Clearly each rational point $a\in L$ yields a purely inseparable map. Let $  V\subset\PP^3$ be the linear span of all  rational
 points $a\in X$ with purely inseparable projection $\pi_a:X\dashrightarrow \PP^2$.
Seeking a contradiction, we assume $L\subsetneqq V$. According to \cite{Kollar 2002}, Lemma 5.1 our  cubic surface $X\subset\PP^3$ can be described
after some change of coordinates by an equation $\lambda y^3+ \sum_{i=1}^3\lambda_iyx_i^2=0$ 
in certain new variables $x_i,y$ for
some scalars $\lambda_i,\lambda\in F$. It follows that the scheme $X$ is reducible, contradiction. This proves the first assertion.

Now suppose that $\O_{X,a}$ is regular, and write  $a=(\alpha_1:\beta_1:\alpha_2:\beta_2)$.
Taking partial derivatives in \eqref{cubic equation}, we see that the tangent plane $T_a(X)\subset\PP^3$ 
is given by the equation $(\alpha_1^2+t_1\beta^2_1)y_1 +(\alpha_2^2+t_2\beta^2_2)y_2 =0$.
Without restriction, we may assume that the second coefficient does not vanish.  In turn, the cubic curve
 $C_a\subset\PP^2$ becomes the zero-locus of 
a polynomial $P(x_1,y_1,x_2)$   divisible by $y_1$. Thus $C_a$ is not integral.
\qed

\medskip
The   $a=(0: 0:1:\lambda)\in X$ with $\lambda\in\FF_4^\times$ show that there are indeed rational points with $\O_{X,a}$ regular 
and $\pi_a:X\dashra\PP^2$ separable.
This reveals that, for regular cubic hypersurfaces in   characteristic two,  the  implication 
$$
\begin{gathered}
\text{$\exists a\in X(F)$ with   $\O_{X,a}$ regular}\\[-.1cm] 
\text{and   $\pi_a:X\dashrightarrow \PP^n$ separable}
\end{gathered}
\quad\Longrightarrow\quad \text{$X$ is unirational}
$$
formulated in the remark on imperfect ground fields in  \cite{Kollar 2002}, page 468,  does not hold without an  additional assumption.
The problem seems to be that all $C_a$ fail to be integral.

\medskip

Let us close the paper with the following observations:
Our cubic surface $X\subset \PP^3$ is Gorenstein, with  $\omega_X=\O_X(-1)$ 
and furthermore $h^0(\O_X)=1$ and $h^1(\O_X)=h^2(\O_X)=0$. In particular, $X$ is a \emph{del Pezzo surface} of degree $K_X^2=3$.

\begin{proposition}
The Picard group $\Pic(X)$ is freely generated by the classes of the invertible sheaves $\O_X(C_1)$ and $\O_X(L)$,
The resulting Gram matrix is
$(\begin{smallmatrix}0	& 2\\2	& -1\end{smallmatrix})$, and  the anticanonical class is given by $-K_X=L+C_1$.
\end{proposition}

\proof  
Let $S\subset \Pic(X)$ be the subgroup generated by the effective Cartier divisors $C_1,L\subset X$.
From the intersection numbers $L^2=-1$, $C_1^2=0$ and $(C_1\cdot L)=2$ we see that $C_1,L\in S$ form a basis,
with the Gram matrix from the assertion. Furthermore, we have 
$-K_X=L+C_1$.  Our task is to show that $S\subset\Pic(X)$ is an equality.

Recall that we have a fibration $f:X\ra \PP^1$. The generic fiber is a twisted form of the double line,
and its Picard group is generated by $\O_{X_\eta}(L)$. 
Likewise, all closed fibers are irreducible, and we conclude that $S \subset\Pic(X)$ has finite index.

Using $\operatorname{disc}(S)=-4$, we see that   the discriminant group  $S^*/S$ has order four.
Write $e_1,e_2\in S$ for the basis  corresponding to the Cartier divisors  $C_1,L\subset X$, and $e_1^*,e_2^*\in S^*$ be the
dual basis. One easily checks that $e_2^*=\frac{1}{2}e_1$ generates the discriminant group.
Seeking a contradiction, we assume that this generator comes from an invertible sheaf  $\shN$.
Then $(\shN\cdot\shN) - (\shN\cdot\omega_X) = \frac{1}{2} (L\cdot C_1) = 1$ is odd. However, this number must be even
by Riemann--Roch, contradiction. Thus $S=\Pic(X)$.
\qed

\medskip
The scheme of non-smoothness $D=\Sing(X/F)$ is disjoint from $L$ and has $\deg(D/\PP^1)$, hence $D$ is linearly equivalent
to $C_1+2L$.
Using this information, we can clarify the occurrence of singularities:

\begin{proposition}
Let $0\leq n\leq 2$ be the dimension of the subvector space generated by the $dt_1,dt_2\in\Omega^1_F$.
Then the scheme $X$ satisfies the regularity condition $(R_n)$, and we have the following implications:
\begin{enumerate}
\item If $n=2$ then the cubic surface $X$ is regular.
\item If $n=1$ then $X$ is normal, and $\O_{X,b}$ is singular for some closed $b\in D$.
\item If $n=0$ then the scheme $X$ is non-normal, with singular locus $\Sing(X)=D$.
\end{enumerate}
\end{proposition}

\proof
Assertion (i) already appeared in Theorem \ref{properties cubic}.
Now suppose that $n=0$, such that both $t_1,t_2\in F$ are squares. After a change of coordinates,
we may assume that $t_1=t_2=1$. Then 
for each rational point of the form $a=(\lambda:\lambda:\mu:\mu)$
the defining polynomial $P=x_1(x_1-y_1)^2 + x_2(x_2-y_2)^2$ lies in the square of the maximal ideal in $\O_{\PP^3,a}$
and it follows that all the local rings $\O_{X,a}$, $a\in D$ are singular. Thus $X$ is singular in codimension one,
hence non-normal. This gives (iii).

Finally, assume that $n=1$. Without restriction, we may assume that $dt_1\neq 0$.
Then $t_1\in F$ is not a square, so the closed point $a=(\sqrt{t_1}:0:1:0)\in X$
is non-rational. Consider the resulting local ring $R=\O_{\PP^3,a}$.
The defining polynomial   \eqref{cubic equation} for the cubic surface obviously lies
in $\maxid_R$ but not in $\maxid_R^2$, hence $\O_{X,a}$ is regular.
If follows that the localization $\O_{X,\zeta}$ is regular as well, where $\zeta$
is the generic point of the scheme of non-smoothness $D=\Sing(X/F)$.
Hence $X$ satisfies $(R_1)$, thus our cubic surface  is normal.

It remains to verify that $\O_{X,b}$ is singular for some closed point $b\in X$.
Suppose for a moment that  $D=\Sing(X/F)$
is non-reduced. Since $D\subset X$ is Cartier, the reduction reduction $E=D_\red$ is another effective Cartier divisor,
and we have $D=nE$ for some integer $n\geq 2$.
However, $D$ is linearly equivalent to $C_1+2E$.
This is primitive in the Picard group, contradiction.
Thus we merely have to check that $D$ is non-reduced.
Its homogeneous coordinate ring is the tensor product $A=A_1\otimes_FA_2$ with factors
$$
A_i=F[x_i,y_i]/(x_i^2-t_iy_i^2),
$$
according to Proposition \ref{locus non-smoothness}.  Consider the field extension $E_1=F(\sqrt{t_1})$.
Then the map $A_1\subset E_1[x_1]$ given by $y_1\mapsto t_1x_1$ is a finite ring extension inside the field of fractions.
Since $dt_2$ is a multiple of $dt_1$, the scalar $t_2\in E_1$ becomes a square,
and we conclude that the rings  $A\subset E_1[x_1,x_2,y_2]/(x_2-\sqrt{t_2}y_2)^2$ are non-reduced.
In turn, the scheme $D=\Proj(A)$ is non-reduced.
\qed

\medskip
Suppose that $t_1,t_2\in F$ are $p$-independent, such that $X$ is a \emph{regular del Pezzo surface
that is not geometrically normal}.
The cone of curves $\operatorname{Eff}(X)$ is the real cone generated by the irreducible curves in the real vector space
$N^1(X)_\RR=\Num(X)\otimes\RR$. In our situation the vector space has rank $\rho=2$, and  contains two  extremal rays, 
which are generated by the fiber $C_1$ and the negative-definite curve $L$,
compare \cite{Kollar 1995}, Lemma 4.12.

In turn there is precisely one minimal model $X\ra Y$, which is  the contraction of $L$.
This is another regular del Pezzo surface that is not geometrically normal.
 Now the degree is $K_X^2=2$, and the anticanonical class
generates   $\Pic(X)=\ZZ$.  
Such examples are interesting, because over fields  of $p$-degree $\pdeg(F)\leq 1$ 
there are no regular del Pezzo surfaces that are not geometrically  normal, 
according to \cite{Fanelli; Schroeer 2020}, Theorem 14.1. 
For more information on del Pezzo surfaces of degree two, we refer to the monographs
of Manin \cite{Manin 1986}, 
Dolgachev \cite{Dolgachev 2012}  and Koll\'ar, Smith and Corti \cite{Kollar; Smith; Corti 2004}.


\end{document}